\documentstyle[amscd,amssymb,verbatim,12pt]{amsart}
\newcommand{\C}{\mathbb{C}}

\newcommand{\R}{\mathbb{R}}

\newcommand{\Aut}{\operatorname{Aut}}
\newcommand{\End}{\operatorname{End}}

\newcommand{\Hom}{\operatorname{Hom}}
\renewcommand{\Im}{\operatorname{Im}}

\newcommand{\Ker}{\operatorname{Ker}}

\newcommand{\SO}{\operatorname{SO}}

\newcommand{\Spin}{\operatorname{Spin}}

\newcommand{\Tr}{\operatorname{Tr}}

\theoremstyle{plain}

\numberwithin{equation}{section}

\begin{document}
\title{The Geometry of Supergravity Torsion Constraints}
\author{John Lott}
\address{Department of Mathematics\\
University of Michigan\\
Ann Arbor, MI  48109-1109\\
USA}
\email{lott@@math.lsa.umich.edu}
\thanks{Research supported by NSF grant DMS-0072154}
\date{August 18, 2001}
\maketitle

\section{Introduction}

This paper, prepared for the 2001 Park City Research Program in
Supergeometry, is an exposition of \cite{Lott (1990)}. References
to relevant
earlier papers can be found in the bibliography of \cite{Lott (1990)}.

As explained in Jim Gates' lectures, an essential ingredient of the
superspace formulation of supergravity is a nonzero
torsion tensor.  Furthermore, there are torsion constraints in the
sense that only certain components of the torsion
tensor are allowed to be nonzero.
These torsion constraints must be stringent enough to give physically
relevant solutions, but flexible enough to allow for nonflat
solutions. For example, for $N = 1$ supergravity on a Lorentzian
3-manifold, Jim Gates motivated the choice that
$T_{\alpha \beta}{}^{c} = i (\gamma^c)_{\alpha \beta}$,
$T_{\alpha b}{}^\gamma$ can be nonzero,
$T_{ab}{}^\gamma$ can be nonzero and the
other torsion components must be zero.

The existence of a nonzero torsion tensor in supergravity theories, and the 
constraints on the torsion tensor, are perhaps surprising from the viewpoint of
standard differential geometry.  In particular, their geometric
meaning is not immediately clear. 

To give a somewhat analogous situation from conventional geometry, 
suppose that
$M^{2n}$ is an almost complex manifold with a Hermitian metric.
Complexifying $TM$ and using notation that will be explained later, 
suppose that we are told that a desirable set of
torsion conditions is
\begin{equation} \label{holo}
T_{ijk}=0, \: \: \: \: T_{ij\bar{k}} = T_{\bar{k}ij}-T_{\bar{k}ji}.
\end{equation}
The
geometric meaning of (\ref{holo}) may also not be immediately clear.
If fact, (\ref{holo}) holds
if and only if $M$ is K\"ahler. The first equation in (\ref{holo})
is the integrability condition
for the almost complex structure and the second equation expresses the
K\"ahler condition, in terms of a unitary basis. Now $M$ is K\"ahler
if and only if near each
$p\in M$, $M$ has the Hermitian geometry of $\mathbb{C}^n$ to first order,
i.e. there exist holomorphic coordinates $\{z^i\}_{i=1}^n$ around $p$
such that the metric tensor takes the form
$g_{i\bar{j}}= \delta_{i\bar{j}}+O(|z|^2)$. Thus
(\ref{holo}) means that $M$ has a first-order flat $U(n)$-geometry at
each point.

We wish to give a similar geometric interpretation for the torsion
constraints of supergravity, as the first-order flatness of
a $G$-structure for some appropriate Lie group $G$. In effect we will do
reverse engineering, taking the known torsion constraints and trying
to find a group $G$ from which they come.
 
The theory of $G$-structures
goes back to \'E. Cartan and was extensively developed in the
1960's.  It is not well known today, perhaps because much of
the literature on $G$-structures is difficult to penetrate.  We will
only discuss the minimal amount of this theory that is needed for
the supergravity torsion constraints.

As our groups $G$ will be super Lie groups, we must first say something
about the structure of super Lie groups.

\section{Super Lie groups}

To start off with an example, what should $GL(p|q)$ mean?
Formally,
\begin{align} \label{gl}
GL(p|q) \text{``$=$''} & \text{ the set of $(p+q) \times (p+q)$ invertible
matrices $\begin{pmatrix} A&B\\C&D\end{pmatrix}$} \\
& \text{with $A$
and  $D$
even, and $B$ and $C$ odd}. \notag
\end{align}
There is the immediate problem that our ground ring is
$\mathbb{R}$ which has no odd elements, so it is not clear how $B$ and
$C$ could be nonzero.  To get around this problem we add
auxiliary odd parameters.
Namely, let $\Lambda$ be a graded-commutative real superalgebra.
Then we define $GL(p|q)(\Lambda)$ as in (\ref{gl}), where
$A$, $B$, $C$ and $D$ now have components in $\Lambda$. This makes perfect
sense, and we can multiply two such matrices to see that
$GL(p|q)(\Lambda)$ is an {\it ordinary}\/ Lie group.

As usual, it is convenient to do 
formal calculations by implicitly thinking that
everything takes value in the unspecified $\Lambda$, but one can also
work in a $\Lambda$-independent setting.  To do so, we look for
a supermanifold, which we will call $GL(p|q)$, with the property that
$GL(p|q)(\Lambda)$ is the set of $\Lambda$-points of $GL(p|q)$.
To construct $GL(p|q)$, consider the affine superspace
${{\R}}^{p^2+q^2|2pq}$. Then we can take
\begin{equation}
GL(p|q)={{\mathbb{R}}}^{p^2+q^2|2pq} \big|_{GL(p)\times GL(q)},
\end{equation}
where we think of $GL(p) \times GL(q)$ as a domain in
$\End({\mathbb{R}}^{p}) \times \End({\mathbb{R}}^{q}) \cong
{\mathbb{R}}^{p^2+q^2}$. As in John Morgan's talk, there is a
group structure on $GL(p|q)$, i.e. a morphism $GL(p|q) \times GL(p|q) 
\to GL(p|q)$, etc. 

To give a concrete description of a general super Lie group $G$ as a
supermanifold, imagine starting with a Lie superalgebra ${\frak g}$, 
the meaning
of which is clear.  Then imagine exponentiating the even subalgebra
${\frak g}^{\text{even}}$. 
One can write out the Jacobi identity for ${\frak g}$
in terms of ${\frak g}^{\text{even}}$ and
${\frak g}^{\text{odd}}$, to obtain four equations. Thinking of these as the
infinitesimal equations for $G$, 
one is led to the following
ingredients for a super Lie group:
\begin{enumerate}
\item An ordinary Lie group $G^{\text{even}}$,
\item A finite-dimensional $G^{\text{even}}$-module $V$ and
\item A $G^{\text{even}}$-equivariant symmetric map 
$d:V \times V \to \mathfrak{g}^{\text{even}}$ such that 
\begin{equation}
d(v_1,v_2)\cdot v_3 + d(v_2,v_3)\cdot v_1 + d(v_3,v_1)\cdot v_2 = 0.
\end{equation}
\end{enumerate}
Given these ingredients, we obtain
a supermanifold $G$ with base space $|G| = G^{\text{even}}$ and
$C^\infty(G)= C^\infty(G^{\text{even}}) \otimes \Lambda^*(V^*)$,
and there is a group structure on $G$. \\ \\
Example : The superlinear group $G=GL(p|q)$ comes from
$G^{\text{even}}=GL(p)\times GL(q)$ and $V=
\Hom({\mathbb{R}}^{p}, {\mathbb{R}}^{q}) \oplus
\Hom({\mathbb{R}}^{q}, {\mathbb{R}}^{p})$. \\ \\
Example :  The superorthogonal group $G=OSp(p|q)$ comes from
$G^{\text{even}}=O(p)\times Sp(q)$ and $V= \Hom({\mathbb{R}}^p,
{\mathbb{R}}^{2q})$. \\ 

Given a supermanifold $M$, we can define a principal $G$-bundle on it
in terms of
transition morphisms
$\phi_{\alpha,\beta} : M \big|_{U_\alpha \cap U_\beta} \to G$, where
$\{U_\alpha\}$ is an appropriate open covering of $|M|$ and
the $\phi_{\alpha,\beta}$'s satisfy a cocycle condition; see John
Morgan's lecture for a similar description of vector bundles.\\ \\
Example : If $M$ is a supermanfold of dimension $(p|q)$ then
the frame bundle $FM \to M$ is a  principal $GL(p|q)$-bundle whose
transition functions are given in terms of local coordinates by
$\phi_{\alpha,\beta}=\left(\frac{\partial z}{\partial 
w}\right)$.\\

Hereafter we will work somewhat formally, but keeping in mind that
super Lie groups should be interpreted as described above.
The first ingredient of a supergeometry is a reduction of $FM$ to some
principal bundle $P \to M$ with some structure group $G$.  That is,
$G$ is a super Lie subgroup of $GL(p|q)$ 
and there is a $G$-equivariant embedding $P \subset FM$.
One can think of $P$ as giving the preferred sets of frames.

For example, if $M$ is an ordinary $n$-dimensional manifold then 
usual Riemannian geometry amounts to a reduction from
the $GL(n)$-bundle $FM$ to an $O(n)$-bundle $OM$. In analogy, if $M$ is a
$(p|2q)$-dimensional supermanifold then one's first attempt to
define a supergeometry might be to take a reduction of $FM$ from a
$GL(p|2q)$-bundle to an $OSp(p|q)$-bundle. This would correspond to having
a superRiemannian metric.  If this were the correct notion of supergeometry
then it would be pretty boring, as it would just be the 
$\mathbb{Z}_2$-grading of what is done in usual Riemannian geometry.  
However, as Jim Gates explained, this is the {\it wrong} notion of
supergeometry, at least from the view of supergravity. Instead,
in the physicists'
description of supergravity, one assumes that $|M|$ is a 
Spin-manifold and one takes $G=\Spin(p)$, together with some
torsion constraints. 

In order to interpret the torsion constraints geometrically 
we will take a slightly different structure group $G$,
but first we must explain the meaning of the torsion tensor.

\section{Torsion}

In Riemannian geometry courses, what we learn about the torsion tensor is
that it's something to be set to zero.  While this is true, it's not very
illuminating.  In a nontrivial sense, 
the real reason that we set the torsion to zero in
Riemannian geometry is because we can.
We now discuss what torsion means in general.  

Let's first consider an ordinary manifold $M^n$. To keep things
concrete, let's take a local basis $\{E^A\}$ of $1$-forms on $M$. 
Following Jim Gates' notation, a $G$-connection can be written as
\begin{equation}
\omega_B{}^A=dz^M \omega_{MB}{}^A = E^C\omega_{CB}{}^A,
\end{equation}
where $\omega_B{}^A$ denotes a $1$-form that takes value in the Lie algebra
${\frak g} \subset
{\frak gl}(n)$. The torsion tensor $T$ is given by
\begin{equation} \label{torsion}
T^A = dE^A + E^B \wedge \omega_B{}^A.
\end{equation} 

As a thought experiment,
given $\{E^A\}$, how much of $T$ can we kill by changing $\omega$? From
(\ref{torsion}), if we send $\omega_B{}^A$ to
$\omega_B{}^A +  \Delta \omega_B{}^A$ then $T$ changes by
\begin{align} \label{kill}
\triangle T^A &= E^B \wedge \triangle \omega_B{}^A\\
&= E^B\wedge E^C \triangle \omega_{CB}{}^A. \notag
\end{align}
The latter is the stuff that we can kill.

To say this more formally, let $W$ be our flat space, so 
$G \subset GL(W)$.
The framing
$\{E^A\}$ at $m$ gives an isomorphism
$T_mM \cong W$. Then we can think of the torsion tensor
$T(m)$ at $m$ as an element of $\Hom(W \wedge W, W)$. Equation (\ref{kill})
defines a linear map
\begin{equation}
\delta : \Hom(W, {\mathfrak{g}}) \to \Hom(W\wedge W, W)
\end{equation}
which sends 
$\triangle \omega$ to $\triangle T$. Note that $\delta$ is defined
purely algebraically. The part of the torsion that we can kill is
$\Im(\delta)$.

Let us define 
\begin{equation}
H^{0,2} = \Hom(W\wedge W,W)/\Im{\delta}.
\end{equation}
In other words, this is the part of the torsion that 
we cannot kill by changing the connection.
Given $T(m) \in \Hom(W\wedge W,W)$, we write its equivalence class in 
$H^{0,2}$ as $[T(m)]$. While we're at it, let's
define ${\frak g}^{(1)}$ to be $\Ker(\delta)$. Given 
$\{E^A\}$ and $T$, from (\ref{torsion})
this is the amount of freedom in the connection $\omega$.

To phrase things in terms of principal bundles, recall that there is the
notion of the soldering form $\tau$, a canonically-defined $W$-valued
$1$-form on $FM$.  
Given the reduction $P \subset FM$, we pullback $\tau$ to $P$ and
give it the same name.  Suppose that 
we have a local section
$s : (U \subset M) \to P$. Then $\{E^A\}$ is just $s^* \tau$.

Let $\omega$ be a connection on $G$, 
i.e. a $G$-equivariant ${\frak g}$-valued
$1$-form on $P$ with the property that $\omega(V_x) \: = \: x$, where
$V_x$ denotes the vector field on $P$ generated by $x \in {\frak g}$. 
The torsion tensor is the horizontal $W$-valued $2$-form on $P$ given by
$T = d \tau + \tau \wedge \omega$.  Using $\tau$, we can also consider
the torsion to be a $G$-equivariant map $T : P \to
\Hom(W\wedge W, W)$. Quotienting by $\Im(\delta)$, we obtain a
$G$-equivariant map $[T] : P \to H^{0,2}$. Note that by construction,
$[T]$ depends only on the reduction $P$, i.e. is independent of the
choice of connection $\omega$.

What is the significance of $[T]$? It gives us an obstruction for $M$ to be
$G$-flat. By $G$-flatness,
 we mean the following. We are given a $G$-reduction of $FM$
to $P$. We assume that the flat space $W$ has a canonical reduction of its 
frame bundle $FW$ to a principal $G$-bundle
$P^{\text{flat}} \subset FW$. Given $m \in M$, consider a
diffeomorphic embedding $\phi : ({\mathcal N} \subset W) 
\to M$, where ${\mathcal N}$ is a
neighborhood of $0 \in W$ and $\phi(0) = m$. We can always
lift $\phi$ to an embedding $\phi_* : F{\mathcal N} \to FM$. 
Here's the geometric question :
does $\phi_*$ send $P^{\text{flat}} \big|_{\mathcal N} \subset
F{\mathcal N}$ to
$P$? If so then for all practical purposes, $M$ is locally the same as $W$.
We say that $M$ is $G$-flat if for each $m \in M$, we can find an
embedding $\phi$ so that $\phi(0) = m$ and $\phi_*$ does send 
$P^{\text{flat}} \big|_{\mathcal N}$ to $P$.

To see what this has to do with the torsion,
let's suppose that $M$ is $G$-flat. Given
$m \in M$, construct $\phi : ({\mathcal N} \subset W) 
\to M$ as above.   
Suppose that $W$ has a $G$-connection $\omega^{\text{flat}}$
with vanishing torsion. (That is, $\omega^{\text{flat}}$ is defined on
$P^{\text{flat}}$.) Then using the embedding $\phi_*$, we
can transfer $\omega^{\text{flat}}$ to obtain a torsion-free connection
defined on $P$ over $\phi({\mathcal N})$. It follows that $[T]$ vanishes in
$H^{0,2}$, at least over $\phi({\mathcal N})$. 

Running the logic backwards,
we see that
a nonvanishing of $[T]$ is an obstruction for $M$ to be $G$-flat. Again,
this is a statement just about the $G$-reduction $P$. In fact, one can define
a notion of $M$ (or more precisely $P$)
being first-order $G$-flat and then a precise 
statement is that
$M$ is first-order $G$-flat at $m \in M$ if and only if $[T(m)]$ vanishes
in $H^{0,2}$ \cite[Theorem 4.1]{Guillemin (1965)}. \\ \\
Example : If $G=O(n) \subset GL(n)$ then one computes algebraically
that $H^{0,2}=0$. Thus there is no obstruction to first-order flatness in
Riemannian geometry. In Lorentzian geometry, this is a form of the
equivalence principle. As ${\frak g}^{(1)} = 0$, there is a unique torsion-free
orthogonal connection, the Levi-Civita connection.\\ \\
Example : If $G=U(n)\subset GL(2n)$ then one finds that $H^{0,2}\ne 0$. 
In fact, the condition for $[T]$ to be zero becomes exactly the equations
in (\ref{holo}). Thus $M$ is first-order $U(n)$-flat if
and only if $M$ is K\"ahler. Again ${\frak g}^{(1)} = 0$, so if $M$ is
K\"ahler then there is a unique torsion-free unitary connection. \\

Now suppose instead that our model space $W$ has a 
constant nonzero torsion $T_0$.
(We will still consider it to be 
a flat space, just one with a nonzero torsion.)
Although less common in conventional geometry than vanishing torsion, this
situation does arise, for example, in CR geometry, the geometry of
hypersurfaces in ${\C}^n$.

If we are to model $M$ by $W$ then we want $M$ to also have a connection
with torsion
$T_0$. This doesn't quite make sense as stated, since 
we still have to take into account the action of $G$. In terms of the local
frame $\{E^A\}$, we want to have a connection whose torsion $T$ 
differs from $T_0$ by a $G$-action, since
then we can perform a gauge transformation to make $T$ identically equal to
$T_0$.
The residual local symmetry is the subgroup $G_0$ of $G$ which preserves
$T_0$. 

In terms of the principal bundle $P$, we have the 
$G$-equivariant map $[T] : P \to H^{0,2}$. Because of the $G$-equivariance
it doesn't make sense to say that $[T]$ lands on $[T_0]$, but it does make
sense to require that for each $p \in P$, $[T](p)$ lies in the $G$-orbit of
$[T_0]$ in $H^{0,2}$. If this is the case then we will say that $M$ (or more
precisely $P$) is first-order $G$-flat. If $M$ is first-order $G$-flat, let
us choose a $G$-connection $\omega$ whose torsion $T : P \to
\Hom(W \wedge W, W)$ takes value in the $G$-orbit of $T_0$. Putting
$P_0 = T^{-1}(T_0)$, we obtain a reduction of $P$ to a principal
$G_0$-bundle $P_0$. In fact, once $T_0$ is nonzero,
it is only natural to make such a reduction.

Finally, let ${\frak g}_0$ denote the Lie algebra of $G_0$ and 
suppose that we have a $G_0$-invariant splitting
${\frak g} = {\frak g}_0 \oplus {\frak g}_1$. When we pullback
$\omega$ to $P_0 \subset P$, it decomposes as
$\omega = \omega^\prime + \omega^{\prime \prime}$.  
Here $\omega^\prime$ is a $G_0$-connection
on $P_0$ and $\omega^{\prime \prime}$ is a tensor.
The torsion equation (\ref{torsion}) becomes
\begin{equation} \label{structure}
T_0^A - E^B \wedge \omega^{\prime \prime}_B{}^A  = dE^A +
E^B \wedge \omega_{B}^{\prime \: A}. 
\end{equation}
That is, although we started with a first-order flat $G$-structure, the
induced $G_0$-structure may not be first-order flat, but instead has
the torsion tensor $T_0^A - E^B \wedge \omega_B^{\prime \prime A}$.

\section{Supergravity torsion constraints}

We return to supergeometry. We first consider unextended, 
i.e. $N = 1$, supergravity theories.

The model flat space $W$ is a superspace ${\R}^{p|q}$ where
${\R}^p$ has an inner product of signature $(p_+, p_-)$ and 
${\R}^q$ is a faithful spinor module for $\Spin(p_+, p_-)$.  
We use the standard notation that lower-case Roman indices are even
indices, Greek indices are odd indices and upper-case Roman indices are
either.
We assume that
there is a charge conjugation operator, i.e. a matrix $C \in
\Aut({\R}^q)$ such that $C \gamma_a C^{-1} = \alpha \gamma_a^T$ and
$C^T = \alpha C$, with $\alpha = \pm 1$. Then as discussed in the
other lectures, to $W$ one can associate a super Poincar\'e group, an 
invariant collection $(D_a, D_\alpha)$ of vector fields and a flat
connection whose only nonzero torsion component is
$(T_0)_{\alpha \beta}^{\: \: \: \: \: c} = (\gamma^c C^{-1})_{\alpha \beta}$. 

We now take $G$ to be a super Lie group of the form
\begin{equation}
G = \left\{ \begin{pmatrix}\rho_1(A) &0\\
\star &\rho_2(A)\end{pmatrix} : A\in \Spin(p_+, p_-), \star \in 
{\mathcal S} \right\}
\end{equation}
where
$\rho_1 : \Spin(p_+, p_-) \to \SO(p_+, p_-)$ is the orthogonal 
representation,  $\rho_2$ is the spinor
representation and $\mathcal{S}$ is 
a $\Spin(p_+, p_-)$-invariant subspace of 
$\Hom({\mathbb{R}}^p,{\mathbb{R}}^q)$. 
That is, $G$ is the semidirect product ${\mathcal S} \widetilde{\times}
\Spin(p_+, p_-)$.
\\ \\
{\bf Claim :} The torsion constraints in supergravity theories of dimension at
most six all arise as the first-order flatness of a $G$-structure.\\

We will work out in detail the three-dimensional
example $(p_+, p_-) = (2, 1)$.  First, we make some general remarks.
There is some freedom in the choice of subspace $\mathcal{S}$. Of course,
because of the $\Spin(p_+, p_-)$-invariance there is only a finite number of
possibilities.
{\it A priori}, different choices of $\mathcal{S}$ can
give different geometries.  In most cases, one can just take 
$\mathcal{S}$ to be all of $\Hom({\mathbb{R}}^p,{\mathbb{R}}^q)$.

The condition of first-order flatness for a $G$-structure
can be written out in equations,
but we will just need the geometric notion.
It is easy to see that the subgroup $G_0$ which preserves $T_0$ is
$\Spin(p_+, p_-)$. In this way we make contact with the physicists'
superspace formulation of supergravity. 
There is an obvious $G_0$-invariant decomposition
${\frak g} = {\frak g}_0 \oplus {\mathcal S}$.
After reducing to a $G_0$-bundle,
the structure equations that
we derive are exactly (\ref{structure}).

The group $G$ preserves the rank-$q$ odd
subspace of ${\R}^{p|q}$. However, if
${\mathcal S} \neq 0$ then it
does not preserve ${\R}^p$.

Geometrically, if we have a first-order flat $G$-structure on 
$M$ then we obtain
\begin{enumerate}
\item A nonintegrable 
odd distribution $D^{\text{odd}}M \subset TM$ of rank $(0|q)$,
\item A spinorial $\Spin(p_+, p_-)$-representation on $D^{\text{odd}}M$,
\item An orthogonal
$\Spin(p_+, p_-)$-representation on $TM/D^{\text{odd}}M$ and
\item A map $[\cdot,\cdot]: D^{\text{odd}}M \times D^{\text{odd}}M
\to TM/D^{\text{odd}}M$, coming from the Lie bracket, which is
conjugate to the corresponding flat space map
${\R}^q \times {\R}^q
\to {\R}^p$ described by $T_0$.
\end{enumerate}

We also endow $TM/D^{\text{odd}}M$ with a compatible 
inner product.
In the case
${\mathcal{S}}
= \Hom({\mathbb{R}}^p,{\mathbb{R}}^q)$, the linear transformations
of $TM$ which preserve the above structure give exactly the group $G$.

The diagonal subgroup $G_0$ of $G$ preserves the subspaces
${\R}^p$ and ${\R}^q$ of ${\R}^{p|q}$. Hence a reduction to a $G_0$-structure
corresponds to a choice of splitting
$TM = (TM/D^{\text{odd}}M) \oplus D^{\text{odd}}M$.

Note that if $\mathcal{S} \neq 0$ then the group $G$ is not a subgroup of 
$OSp \left( p_+, p_- | \frac{q}{2} \right)$.  
This shows again that the notion of a superRiemannian metric
is irrelevant for supergravity theories.

To summarize, to form a supergeometry, suppose that we are given 
${\mathcal S}$.
Then
\begin{enumerate}
\item Pick a reduction $P$ of $FM$ to a $G$-structure, i.e. a set of
framings $\{E^A\}$.
\item Check whether the reduction is first-order flat.  If it isn't,
throw it out.
\item Choose a $G$-connection with the correct torsion.
\item Reduce to the subgroup $G_0$.
\item Write out the structure equations.  Analyze their consistency using
the Bianchi identities.
\end{enumerate}

The condition of being first-order flat in 2 involves first-derivatives of
the frame $\{E^A\}$. It is not at all obvious how to parametrize the set of
solutions. One wants to do so in order to find the
independent supergravity fields. In the case of three-dimensional
supergravity, Jim Gates explained how one can parametrize the
independent fields using the
spinorial frame $\{E^\alpha\}$. In the four-dimensional case, one
encounters prepotentials $\{H^A\}$. It appears that one must do a 
case-by-case analysis to find the independent fields.

There is generally not a unique choice of connection in 3, as 
${\frak g}^{(1)} \neq 0$. However, the ambiguity is mild in the sense that
different choices of connection lead to the same structure equations,
as we will see in the three-dimensional case. 

There is a strong analogy between supergravity theory and CR manifolds.
In fact, there is a dictionary
\begin{align}
W & \longleftrightarrow (S^{2N-1} \subset {\C}^N) \\
D^{\text{odd}}M & \longleftrightarrow D^{\text{complex}}M \notag \\
T_0 & \longleftrightarrow \text{ the Levi form of 
$S^{2N-1} \subset {\C}^N$} \notag \\
\text{ superconformal geometry } & \longleftrightarrow 
\text{ CR geometry (Chern-Moser) \cite{Chern-Moser (1975)}} \notag \\
\text{ supergravity } & \longleftrightarrow 
\text{ pseudoHermitian geometry (Webster)
\cite{Webster (1978)}} \notag
\end{align}
Of course, in CR geometry one does not have odd variables, but the
role of the odd distribution $D^{\text{odd}}M$ is played by the complex 
distribution
$D^{\text{complex}}M \subset TM$.
It is a historical coincidence that Chern and Moser were
analyzing CR manifolds around the same time and in almost the same way
that Wess and Zumino were deriving the superspace formulation of supergravity
theories \cite{Grimm-Wess-Zumino (1979),Wess-Zumino (1977)}.
The second half of the Chern-Moser paper looks at the
structure equations of a CR manifold, 
chooses preferred connections and analyzes the
consequences of the Bianchi identities in a way that mirrors the Wess-Zumino
work, although in a very different language.

To deal with the case of extended supergeometries, let $K$ be a 
Lie group. We assume that ${\R}^q$ is the tensor product of
representation spaces of $\Spin(p_+, p_-)$ and $K$. We take 
${\mathcal S}$ to be a $(\Spin(p_+, p_-) \times K)$-invariant subspace
of $\Hom({\R}^p, {\R}^q)$. Then we put $G = {\mathcal S} \widetilde{\times}
(\Spin(p_+, p_-) \times K)$ and proceed as before. Note that the group
$K$ is gauged.

Our claim is that the torsion constraints of all supergravity theories of
dimension at most six arise from the above procedure.  To be precise,
this is true for theories with an offshell superspace formulation, i.e.
for $q \: \le \: 8$. Also,
in four dimensions there are various superspace formulations of the
same onshell action and we only pick up the ``minimal'' superspace
formulation.

The claim is shown explicitly in \cite{Lott (1990)} for the following 
supergravity theories :

\begin{equation}
\begin{matrix}
\underline{p} & \underline{p_+} & \underline{p_-} & \underline{q} & 
\underline{N} \\
& & & & \\
1 & 1 & 0 & 1 & 1 \\
2 & 1 & 1 & 1 & (1,0) \\
2 & 1 & 1 & 2 & (1,1) \\
2 & 1 & 1 & 2 & (2,0) \\
2 & 1 & 1 & 3 & (2,1) \\
2 & 2 & 0 & 4 & (2,2) \\
3 & 2 & 1 & 2 & 1 \\
4 & 3 & 1 & 4 & 1 \\
4 & 4 & 0 & 8 & 2 \\
6 & 5 & 1 & 8 & 1
\end{matrix}
\end{equation}

We now work out the three-dimensional case in some detail.

\section{Three-dimensional supergravity}

We take $p_+ = 2$, $p_- = 1$ and $q = 2$. Then $\Spin(p_+, p_-) =
SL(2, {\R})$ and ${\R}^q$ has the standard $SL(2, {\R})$-representation.  

We use a notation in which an element
$M \in sl(2, {\R})$ is represented by a traceless matrix $M_\pm^{\: \: \: \pm}$.
Using the $SL(2, {\R})$-invariant 
symplectic form $\epsilon$ on ${\R}^2$ to raise and lower indices,
we can represent $M$ as a symmetric matrix $M_{\pm \pm}$. We also identify
$S^2({\R}^2)$ with Minkowski $3$-space, to write an element $P \in {\R}^3$ as
a symmetric matrix $P^{\pm \pm}$.

Let us first take ${\mathcal S} = \Hom({\R}^3, {\R}^2)$. 
One has $\dim(V) = 5$, $\dim({\frak g}) = 9$,
$\dim(\Hom(V, {\frak g})) = 45$, 
$\dim(\Hom(V \wedge V, V)) = 60$,
$\dim( {\frak g}^{(1)}) = 12$, $\dim(\Im(\delta)) = 33$ and 
$\dim(H^{0,2}) = 27$. It turns out that $G$ acts trivially on 
$H^{0,2}$. 

Suppose that we have a $G$-structure. For it to be first-order flat, $[T]$ must
equals $[T_0]$ identically.  Suppose that this is the case. Then we
choose a connection $\omega$ and reduce to a $G_0$-structure,
where $G_0 = SL(2, {\R})$.

The structure equations are given by  (\ref{structure}). To write these
explicitly, we write the $sl(2,{\R})$-valued connection 
$1$-form $\omega^\prime$ in its spinor
representation as $\omega_\beta^{\prime \alpha}$ and we write it in its
orthogonal representation as
$\omega_{\epsilon \phi}^{\prime \alpha \beta} =
\omega_{\epsilon}^{\prime \alpha} \delta_\phi^{\: \: \beta} + 
\omega_{\phi}^{\prime \alpha} \delta_\epsilon^{\: \: \beta}$.
Then (\ref{structure}) becomes
\begin{align} \label{structure2}
 - E^\alpha \wedge E^\beta & = 
dE^{\alpha \beta} +  E^{\gamma \beta} \wedge 
\omega_{\gamma}^{\prime \: \alpha} + E^{\alpha \gamma} \wedge
\omega_{\gamma}^{\prime \: \beta}, \\
- E^{\gamma \delta} \wedge \omega_{\delta \gamma}^{\prime \prime \: \alpha}  
& = 
dE^{\alpha} +  E^{\beta} \wedge 
\omega_\beta^{\prime \: \alpha}. \notag
\end{align} 

Let us write the $1$-form $\omega_{\delta \gamma}^{\prime \prime \: \alpha}$
as
\begin{equation} \label{1form}
\omega_{\delta \gamma}^{\prime \prime \: \alpha} =
- E^{\epsilon} \: T_{\epsilon, \delta \gamma}^{\: \: \: \: \: \: \:
\alpha}
- \frac12
E^{\epsilon \phi} \: T_{\phi \epsilon, \delta \gamma}^{\: \: \: \: \: \:
\: \: \: \alpha}.
\end{equation}
Then (\ref{structure2}) becomes
\begin{align} \label{structure3}
dE^{\alpha \beta} +  E^{\gamma \beta} \wedge 
\omega_{\gamma}^{\prime \: \alpha} + E^{\alpha \gamma} \wedge
\omega_{\gamma}^{\prime \: \beta} & =
 - E^\alpha \wedge E^\beta, \\
dE^{\alpha} +  E^{\beta} \wedge 
\omega_\beta^{\prime \: \alpha} & = 
E^{\gamma \delta} \wedge E^\epsilon \:
T_{\epsilon, \delta \gamma}^{\: \: \: \: \: \: \:
\alpha}
+ \frac12 E^{\gamma \delta} \wedge
E^{\epsilon \phi} \: T_{\phi \epsilon, \delta \gamma}^{\: \: \: \: \: \:
\: \: \: \alpha}. \notag
\end{align} 

We recognize these as 
the structure equations for three-dimensional supergravity
with structure group $SL(2, {\R})$ and nonzero torsion components
$T_{\epsilon, \delta}^{\: \: \: \: \:
\alpha \beta}$,
$T_{\epsilon, \delta \gamma}^{\: \: \: \: \: \: \:
\alpha}$ and
$T_{\phi \epsilon, \delta \gamma}^{\: \: \: \: \: \:
\: \: \: \alpha}$, as desired.

Recall that there was an ambiguity of ${\frak g}^{(1)}$ in the choice of
the connection.  One can check that this amounts to changing
$T_{\phi \epsilon, \delta \gamma}^{\: \: \: \: \: \:
\: \: \: \alpha}$ in (\ref{1form}) by something of the form
$U_{\phi \epsilon, \delta \gamma}^{\: \: \: \: \: \:
\: \: \: \alpha}$ with 
$U_{\phi \epsilon, \delta \gamma}^{\: \: \: \: \: \:
\: \: \: \alpha} = 
U_{ \delta \gamma, \phi \epsilon}^{\: \: \: \: \: \:
\: \: \: \alpha}$. However, as
$E^{\gamma \delta} \wedge
E^{\epsilon \phi} = -  
E^{\epsilon \phi} \wedge E^{\gamma \delta}$, the structure equation
(\ref{structure3}) remains unchanged.

The Bianchi identities imply that one can express the torsion and curvature
in terms of a function $R$ and a tensor $G_{\alpha \beta \gamma}$ which is
totally symmetric in its indices.  A calculation gives
\begin{align}
T_{\epsilon, \delta \gamma}^{\: \: \: \: \: \: \:
\alpha}  = 
R ( & \epsilon_{\epsilon \delta} \delta_\gamma^{\: \: \alpha} +
\epsilon_{\epsilon \gamma} \delta_\epsilon^{\: \: \alpha}), \\
T_{\phi \epsilon, \delta \gamma}^{\: \: \: \: \: \:
\: \: \: \alpha} =
\frac12 ( &
\epsilon_{\phi \delta} \delta_\epsilon^{\: \: \alpha} \nabla_\gamma R +
\epsilon_{\epsilon \delta} \delta_\phi^{\: \: \alpha} \nabla_\gamma R +
\epsilon_{\phi \gamma} \delta_\epsilon^{\: \: \alpha} \nabla_\delta R +
\epsilon_{\epsilon \gamma} \delta_\phi^{\: \: \alpha} \nabla_\delta R \notag \\
& -
\epsilon_{\delta \phi} \delta_\gamma^{\: \: \alpha} \nabla_\epsilon R -
\epsilon_{\gamma \phi} \delta_\delta^{\: \: \alpha} \nabla_\epsilon R -
\epsilon_{\delta \epsilon} \delta_\gamma^{\: \: \alpha} \nabla_\phi R -
\epsilon_{\gamma \epsilon} \delta_\delta^{\: \: \alpha} \nabla_\phi R \notag \\
& +
G_{\phi \delta}^{\: \: \: \: \alpha} \epsilon_{\epsilon \gamma} +
G_{\epsilon \delta}^{\: \: \: \: \alpha} \epsilon_{\phi \gamma} +
G_{\phi \gamma}^{\: \: \: \: \alpha} \epsilon_{\epsilon \delta} +
G_{\epsilon \gamma}^{\: \: \: \: \alpha} \epsilon_{\phi \delta}
), \notag
\end{align}
with the constraint
\begin{equation}
\nabla_\alpha G_{\beta \gamma}^{\: \: \: \: \alpha} +
\nabla_\beta \nabla_\gamma R + \nabla_\gamma \nabla_\beta R = 0. 
\end{equation}
(We use different conventions than Jim Gates, but the results are
equivalent.) 

Now suppose that we instead take ${\mathcal S}$ to be the subspace of
$\Hom({\R}^3, {\R}^2)$ consisting of maps $M$ that can be written in the form
$M_{\delta \gamma}^{\: \: \: \: \alpha} = 
Z_\delta \delta_\gamma^{\: \: \alpha} + 
Z_\gamma \delta_\delta^{\: \: \alpha}$ for some $Z$. That is,
${\mathcal S}$ consists of the maps $M \in \Hom(S^2({\R}^2), {\R}^2)$ with the
property that there exists a
$z \in ({\R}^2)^*$ such that $M(v,w) = z(v) w + z(w) v$.
Then it turns out that we obtain the same geometry as if we had
taken ${\mathcal S}$ to be all of $\Hom({\R}^3, {\R}^2)$ 
\cite[Proposition 14]{Lott (1990)}.
Next, suppose that we take
${\mathcal S}$ to be the subspace of
$\Hom({\R}^3, {\R}^2)$ consisting of maps $M$ such that
$M_{\alpha \beta}^{\: \: \: \: \beta} = 0$. That is, if we define
$M_v : {\R}^2 \rightarrow {\R}^2$ by
$M_v(w) = M(v \otimes w + w \otimes v)$ then
${\mathcal S}$ consists of the maps $M \in \Hom(S^2({\R}^2), {\R}^2)$ such that
for all $v \in {\R}^2$, $\Tr(M_v) = 0$. In this case it turns out
that the geometry we obtain is equivalent to that obtained from
taking ${\mathcal S}$ to be all of $\Hom({\R}^3, {\R}^2)$, but setting
the superfield $R$ to be zero
\cite[Proposition 13]{Lott (1990)}.
Finally, if ${\mathcal S} = 0$ then we only obtain flat solutions.

\section{Further topics}

\subsection{Higher order obstructions to flatness}

In the theory of $G$-structures, it is generally not true that
first-order flatness implies flatness.  It is true when
$G = GL(n, {\C}) \subset GL(2n, {\R})$, as the first-order flatness
amounts to the vanishing of the Nijenhuis tensor and this implies the
integrability of the complex structure.  It is not true when
$G = O(n) \subset GL(n, {\R})$, as the Riemann curvature tensor is an
obstruction to flatness.

There is a algebraic theory of higher order obstructions to flatness, which
live in the so-called 
Spencer cohomology groups $H^{i,2}$. A clear exposition of this
theory is in \cite{Guillemin (1965)}. For example, if $G = O(n)$ then
$H^{i,2}$ vanishes if $i \neq 1$, while
$H^{1,2}$ consists of the tensors with the symmetry of the Riemann curvature
tensor.  One of the main issues in the theory of $G$-structures is to know
when the vanishing of all of the algebraic obstructions to flatness
actually implies flatness. For example, if $G = GL(n, {\C}) \subset GL(2n, {\R})$
then the flatness is the Newlander-Nirenberg theorem.

In the case of the supergeometry group $G$ with ${\mathcal S} = 
\Hom({\R}^p, {\R}^q)$, it turns out that if $i > 0$ then the only nonvanishing
Spencer cohomology group is $H^{1,2}$, which also consists of the 
tensors with the symmetry of the Riemann curvature
tensor for a $p$-manifold \cite[Propositions 18,19]{Lott (1990)}.  In other
words, there are no formal obstructions to flatness beyond the
curvature tensor $R_{abcd}$ with even indices.

\subsection{Superconformal geometry}

In Riemannian geometry, the Weyl tensor fits nicely into
the framework of Cartan connections.  The latter means that one has Lie groups
$H \subset G$, a principal $H$-bundle $P \rightarrow M$ and a
${\frak g}$-valued $1$-form $\omega$ on $P$ such that \\
1. $\omega$ is $H$-equivariant, \\
2. For all $x \in {\frak h}$, 
$\omega(V_x) = x$, where $V_x$ is the vector field 
on $P$ generated by $x$, and \\
3. For all $p \in P$, $\omega$ gives an isomorphism from 
$T_p P$ to ${\frak g}$.

Define the curvature $\Omega$ as usual to be $d \omega + \omega^2$.
Suppose that $\Omega = 0$. For simplicity, we assume that $P$ and $H$ are
connected. If $p_0 \in P$ 
is a basepoint then we take the path-ordered integral 
of $\omega$ along paths from $p_0$. This gives a map from the universal
cover $\widetilde{P}$ to $G$. Taking a quotient, we obtain a 
$\pi_1(M)$-equivariant map $\alpha$ from
$\widetilde{M}$ to $G/H$.
By condition 3, $\alpha$ is a local diffeomorphism. Thus $\alpha$ is a
developing map and we have coordinate charts on $M$ modeled by
domains in $G/H$, with each transition map
coming from the left-action of an element of $G$.

Of course, in general we cannot assume that $\omega$ is flat.
One can write Riemannian geometry in terms of Cartan connections by
taking $G$ to be the Euclidean group 
${\R}^n \widetilde{\times} SO(n)$ and $H = SO(n)$. Then
the ${\R}^n$-component of $\omega$ can be identified with the soldering form
and the $so(n)$-component of $\omega$ can be identified with the Riemannian
connection.  By an appropriate choice of the Riemannian connection,
i.e. choosing
the Levi-Civita connection, we can kill the ${\R}^n$-component of $\Omega$,
i.e. the torsion.
The remaining $so(n)$-component of $\Omega$ is the Riemannian curvature, which
we see as an obstruction to Euclidean flatness.

In the case of $n$-dimensional conformal geometry, $G$ is the conformal group
$SO(n+1, 1)$. It acts on $S^n$ by conformal transformations. Fix a point
$\infty \in S^n$ and let $H$ be the stabilizer of $\infty$.  Equivalently,
writing ${\R}^n = S^n - \infty$, $H$ is the conformal group of ${\R}^n$.
By construction, $G/H = S^n$. 

Algebraically, ${\frak g}$ is a graded Lie algebra ${\frak g}^{(-1)} \oplus
{\frak g}^{(0)} \oplus {\frak g}^{(1)}$, with ${\frak g}^{(-1)} \cong
{\R}^n$, ${\frak g}^{(0)} \cong o(n) \oplus {\R}$ and ${\frak g}^{(1)} \cong
{\R}^n$. Then ${\frak h} = {\frak g}^{(0)} \oplus {\frak g}^{(1)}$.
Given a Cartan connection $\omega$, we decompose it as
$\omega = \omega^{(-1)} + \omega^{(0)} + \omega^{(1)}$, and similarly for
its curvature $\Omega$. We
can identify the ${\R}^n$-valued $1$-form
${\omega}^{(-1)}$ with the soldering form for $M$. 

To see the Weyl curvature as an obstruction to conformal flatness, one
assumes that one is given the soldering form, i.e.
${\omega}^{(-1)}$, and the Levi-Civita connection, i.e. the
component of $\omega^{(0)}$ in $o(n)$. 
Then the question is whether one can extend
these components to form a Cartan connection $\omega$ which is flat. Taking
the ${\R}$-component of $\omega^{(0)}$ 
to vanish, one can use the freedom in
${\omega}^{(1)}$ to make $\Omega^{(0)}$
equal to the Weyl tensor. In this
way, one sees that the Weyl curvature is an obstruction to conformal flatness
of a Riemannian metric.  If $n = 3$ then the Weyl curvature vanishes but
${\Omega}^{(1)}$ gives its three-dimensional
analog.

It is of interest to treat superconformal geometry in terms of Cartan
connections.  It is a remarkable fact that
the superconformal groups $G$ are simple super Lie
groups, with graded Lie superalgebra
${\frak g}^{(-1)} \oplus {\frak g}^{(-1/2)} \oplus
{\frak g}^{(0)} \oplus {\frak g}^{(1/2)} \oplus {\frak g}^{(1)}$. Here
${\frak g}^{(-1)} \cong {\R}^p$, ${\frak g}^{(-1/2)} \cong {\R}^q$,
${\frak g}^{(0)} \cong o(p_+, p_-) \oplus {\R} \oplus k$, 
${\frak g}^{(1/2)} \cong {\R}^q$ and ${\frak g}^{(1)} \cong {\R}^p$, where
$k$ is the Lie algebra of an internal symmetry group. 
We take $H$ to be
the subgroup with Lie algebra 
${\frak g}^{(0)} \oplus {\frak g}^{(1/2)} \oplus {\frak g}^{(1)}$. 

We decompose a Cartan connection $\omega$ as
\begin{equation}
\omega = \omega^{(-1)} + \omega^{(-1/2)} + \omega^{(0)} + \omega^{(1/2)}
+ \omega^{(1)},
\end{equation}
and similarly for its curvature $\Omega$.  
Now ${\omega}^{(-1)} + {\omega}^{(-1/2)}$
can be identified with the soldering form $\{E^a, E^\alpha\}$. Suppose
that we are given ${\omega}^{(0)}$
and that ${\Omega}^{(-1)}$ vanishes.
(Note that $\Omega^{(-1)}$
only depends on ${\omega}^{(-1)}$, ${\omega}^{(-1/2)}$ and
${\omega}^{(0)}$. Its vanishing corresponds to having the flat space
expression for the torsion component $T^a$.) 
Then the question is whether we can
extend these components to a Cartan connection with vanishing curvature.
In general one cannot, but one can choose $\omega^{(1/2)}$ and
$\omega^{(1)}$ so that certain components of $\Omega$ vanish.  The
remaining components are the obstruction to conformal flatness.

For example, in the case of three dimensions, $G = OSp(1|2)$. It turns out
that one can uniquely choose $\omega^{(1/2)}$ and
$\omega^{(1)}$ so that $\Omega^{(- 1/2)} = \Omega^{(0)} = 0$
\cite[Proposition 23]{Lott (1990)}. The remaining curvature components,
$\Omega^{(1/2)}$ and $\Omega^{(1)}$, are the supersymmetric analog of the
three-dimensional conformal tensor.

\subsection{Onshell theories}

For supergravity theories with $N$ supersymmetries,
if $N$ is large enough then it turns out that the superspace
torsion constraints already imply the equations of motion, i.e. that the
theory is onshell.  This is the case for four-dimensional supergravity
if $N > 2$. In these cases the torsion constraints do not follow the pattern
that we have described above, 
and we do not know of their geometric interpretation.

\subsection{SuperK\"ahler geometry}

SuperK\"ahler manifolds $M$ whose base space
$|M|$ has one complex dimension are
well understood.  As a purely mathematical question, one can ask about
the higher-dimensional situation.  Of course, we are not thinking of the
approach of defining
superK\"ahler forms on a supercomplex manifold, but rather of
applying the theory
of $G$-structures. 

If $|M|$ has $n$ complex dimensions then
there is a natural spinor representation
$\rho_2: U(1) \times SU(n) \rightarrow \End( \Lambda^{*,0} ({\C}^n))$ of
real dimension $2^{n+1}$.
One's first attempt to define a superK\"ahler geometry might be to
require first-order flatness of an ${\mathcal S} \widetilde{\times}
(U(1) \times SU(n))$-structure, where ${\mathcal S}$ is a
$(U(1) \times SU(n))$-invariant subspace of
$\Hom \left( {\R}^{2n}, {\R}^{2^{n+1}} \right)$. However, one finds that
this gives a flat geometry even in the case when $|M|$ has
one complex dimension
\cite[Proposition 27]{Lott (1990)}. Instead, it turns out that one must use
the additional ``chiral'' action of ${\C}^*$ on $\Lambda^{*,0} ({\C}^n)$ which
multiplies an even form by $z \in {\C}^*$ and multiplies an odd form by
$\overline{z}^{-1}$.

Thus we try taking $G = {\mathcal S} \widetilde{\times}
({\C}^* \times U(1) \times SU(n))$. If $|M|$ has one complex dimension then
one finds that this gives the right answer. In fact, one obtains the
same geometry whether one takes 
${\mathcal S} = \Hom \left( {\R}^{2}, {\R}^{4} \right)$ or
${\mathcal S} = \Hom_{{\C}} \left( {\C}, {\C}^{2} \right)$
\cite[Proposition 26]{Lott (1990)}.

If $|M|$ has complex dimension two then one finds that having
a first-order flat $G$-structure, with $G = 
\Hom \left( {\R}^{4}, {\R}^{8} \right) \widetilde{\times}
({\C}^* \times U(1) \times SU(2))$, implies that $|M|$ is a Hermitian
locally symmetric space \cite[Proposition 28]{Lott (1990)}.
(Note that $G$ is a subgroup of the structure group of a
four-dimensional $N=2$ Riemannian
supergeometry.) It may be that this is the best 
that one can do.  In \cite[Proposition 30]{Lott (1990)} we explored
the consequences of relaxing the torsion constraints, but the results were
inconclusive.  In any event, we do not have a general understanding
of superK\"ahler geometry.

We take this opportunity to correct some mistakes in \cite{Lott (1990)}.
Throughout \cite{Lott (1990)} we wrote $\End({\R}^p, {\R}^q)$ when we should
have written $\Hom ({\R}^p, {\R}^q)$. Equation (146) of \cite{Lott (1990)}
should read $\omega_{\Theta_1}^{\: \: \: \Theta_1} + 
\omega_{\Theta_2}^{\: \: \: \Theta_2} = 
\omega_{z}^{\: \: z}$.

I thank the participants of the Park City Research Program for discussions.
I especially thank Dave Morrison for providing his TEX notes of the
Park City talk.


\begin{thebibliography}{100}

\ifx\undefined\bysame
\newcommand{\bysame}{\leavevmode\hbox to3em{\hrulefill}\,}
\fi
\bibitem{Chern-Moser (1975)} S. Chern and J. Moser,
{\em Real hypersurfaces in complex manifolds},
Acta Math. {\bf 133} (1975), 219--271.

\bibitem{Grimm-Wess-Zumino (1979)} R. Grimm, J. Wess and B. Zumino,
{\em A complete solution of the Bianchi identities in superspace with 
supergravity constraints},
Nuclear Phys. B{\bf 152} (1979), 255--265.

\bibitem{Guillemin (1965)} V. Guillemin, 
{\em The integrability problem for $G$-structures},
Trans. Amer. Math. Soc. {\bf 116} (1965), 544--560.

\bibitem{Lott (1990)} J. Lott, 
{\em Torsion constraints in supergeometry},
  Comm. Math. Phys. {\bf 133} (1990), 563--615.

\bibitem{Webster (1978)} S. Webster,
{\em Pseudo-Hermitian structures on a real hypersurface},
J. Diff. Geom. {\bf 13} (1978), 25--41.

\bibitem{Wess-Zumino (1977)} J. Wess and B. Zumino,
{\em Superspace formulation of supergravity},
Phys. Lett. B{\bf 66} (1977), 361--364.
\end{thebibliography}
\end{document}